\documentclass{amsart}

\usepackage{amssymb}
\usepackage{amsmath}
\usepackage{array}
\usepackage{graphicx}
\usepackage{caption}
\usepackage{subcaption}

\newcommand\liegr{\sf}

\newcommand{\SU}[1]{\mbox{${\liegr SU}(#1)$}}

\newcommand{\U}[1]{\mbox{${\liegr U}(#1)$}}
\newcommand{\SP}[1]{\mbox{${\liegr Sp}(#1)$}}
\newcommand{\SO}[1]{\mbox{${\liegr SO}(#1)$}}

\newcommand{\OG}[1]{\mbox{${\liegr O}(#1)$}}
\newcommand{\Spin}[1]{\mbox{${\liegr Spin}(#1)$}}

\newcommand\fieldsetc{\mathbb}

\newcommand{\Z}{\fieldsetc{Z}}
\newcommand{\R}{\fieldsetc{R}}
\newcommand{\C}{\fieldsetc{C}}
\newcommand{\Q}{\fieldsetc{H}}

\newtheorem{thm}{Theorem}[section]

\theoremstyle{remark}

\title{Representations of low copolarity}

\author[A. M. de S\'a Gomes]{Andr\'e Magalh\~aes de S\'a Gomes}

\author[C. Gorodski]{Claudio Gorodski}

\address{Instituto de Matem\'atica e Estat\'\i stica, Universidade de
S\~ao Paulo, Rua do Mat\~ao, 1010, S\~ao Paulo, SP 05508-090, Brazil}

\email{andremsg@usp.br}

\address{Instituto de Matem\'atica e Estat\'\i stica, Universidade de
S\~ao Paulo, Rua do Mat\~ao, 1010, S\~ao Paulo, SP 05508-090, Brazil}

\email{gorodski@ime.usp.br}

\thanks{The first author has been supported by the CAPES
  scholarship 88882.377939/2019-01, and the second author has been partially
  supported by the CNPq grant 302882/2017-0
  and the FAPESP project 2016/23746-6.}

\date{\today}

\begin{document}

\begin{abstract}
  We classify irreducible representations of compact connected
  Lie groups whose orbit space is isometric to the orbit space of a
  representation of a compact Lie group of dimension~$7$, $8$ or $9$.
  They turn out to be closely related to symmetric spaces, with one
  exception only. 
\end{abstract}

\maketitle

\section{Introduction}

A representation $\rho:G\to \OG V$ of a compact Lie group~$G$
is called \emph{polar} if there
exists a subspace $\Sigma$, called a \emph{section},
that meets all $G$-orbits, and meets them
always orthogonally. This notion was made explicit over 35 years ago
and has many connections and ramifications, among them, with symmetric spaces,
invariant theory and submanifold geometry (see e.g.~\cite{D,DK,BCO}).
For a polar representation
as above, the space of orbits $V/G$ can be recovered, as a metric space,
as the space of orbits of a finite group action, namely,
$\Sigma/\mathcal W$, where $\mathcal W$
is the largest subquotient group of $G$ acting on the section $\Sigma$.
Indeed, as an easy consequence of O'Neill's formula for Riemannian
submersions, one sees that the later property characterizes polar
representations (cf.~\cite{GL2}, introd.). 

Polar representations were classified by Dadok~\cite{D};
it follows from this classification that a polar representation of a
connected group has the same orbits as the isotropy representation
of a symmetric space. 
The above considerations and the general importance of polar representations
led the authors of~\cite{GL} to seek other classes of orthogonal representations
whose orbit spaces can be presented in different ways, and to investigate the
mysterious interplay between geometric and algebraic aspects of
representations. 
Namely, one says that two representations are \emph{quotient-equivalent} if they have isometric orbit spaces; in addition, 
if the underlying group of one representation has dimension strictly
less than the dimension of the underlying group of the other one, then one says
the former representation is a \emph{reduction} of the latter one. It follows
that the polar representations of connected groups are precisely those
that admit a reduction to a finite group action.

A \emph{minimal reduction} of a given representation is a reduction whose
underlying group has minimal dimension in that quotient-equivalence class;
in this case, the dimension of this group is called
the \emph{abstract copolarity} of the given representation. For instance,
polar representations of connected groups are precisely those
representations of abstract copolarity zero. 
As one of the main results of~\cite{GL}, it was shown that 
if a non-polar irreducible representation of a compact connected Lie group
admitting a reduction has abstract copolarity at most six, then
the representation is \emph{toric}, that is,
it has a minimal reduction to a representation of a finite
extension of a torus. 
In the same paper, it was found a counterexample for
abstract copolarity $7$: $(\U3\times\SP2,\C^3\otimes_{\mathbb C}\C^4)$
reduces to a $\Z_2$-extension of (the $7$-dimensional group)
$(\OG3\times\U2,\R^3\otimes_{\mathbb R}\R^4)$,
where $\OG3\times\U2$ sits inside $\U3\times\SP2$ as a symmetric
subgroup (see also~\cite{Gom}). 
Toric irreducible representations were later classified
in~\cite{GL2} (they are mostly related to Hermitean symmetric spaces;
see also~\cite{Pa}). 

In the present work, we wanted to understand the extent of the above
counterexample. As it turns out, we can show it is the only counterexample
in abstract copolarity $7$, and we go a little further:

\begin{thm}\label{main}
  Let $\rho:G\to \OG V$ be a non-polar, non-reduced,
  irreducible representation of a compact connected Lie group $G$
  on the Euclidean space $V$. Assume the abstract copolarity of~$\rho$
  is $7$, $8$ or $9$. Then $\rho$ is either toric, quaternion-toric,
  or equivalent to~$(\U3\times\SP2,\C^3\otimes_{\mathbb C}\C^4)$.
\end{thm}

A representation is called \emph{non-reduced} if it admits a reduction;
it is called \emph{quaternion-toric} (or \emph{q-toric}, for short)
if it is non-polar and admits a reduction to a
representation of a group whose identity component
is isomorphic to~$\SP1^k$ for some $k>0$. Q-toric representations
were classified in the irreducible case in~\cite{GG} (they are
related to quaternion-K\"ahler symmetric spaces); in particular,
it was shown that $k=3$ always holds, so they have abstract copolarity $9$.
In particular, an irreducible representation of a connected
Lie group with abstract copolarity~$8$ must be toric. 

Another consequence of Theorem~\ref{main} is that, for the
representations in the theorem, the abstract copolarity
coincides with the copolarity. The \emph{copolarity} of a 
representation $\rho:G\to \OG V$ is the minimal possible
dimension of~$H$, where $\tau:H\to \OG W$ is an arbitrary
reduction of $\rho$ such that
$H$ is a subgroup of $G$, $W$ is a subspace of $V$, and $\tau$ is the
restriction of $\rho$. Thus the copolarity considers only
reductions ``embedded'' in the given representation; in particular,
the abstract copolarity is bounded above by the copolarity. 
It is an open problem to decide whether they always coincide. 

In view of the results and calculations in the present paper, one is tempted
to formulate the following conjecture: \emph{Suppose $G$ is connected and 
  $\rho:G\to \OG V$ is a non-polar, non-reduced,
  irreducible representation that reduces to $\tau:H\to \OG W$.
  Then $H^0$ is not simple.}

The results in this paper are part of the PhD thesis of the first author.

\section{Preliminaries}\label{prelim}

Let $\rho:G\to \OG V$ be a representation of a compact Lie group $G$.
Then $V$ inherits a $G$-invariant stratification by~\emph{orbit types},
namely, two points in $V$ are in the same stratum if and only if they have
conjugate isotropy groups. This stratification projects to a stratification
of the orbit space $X=V/G$. The strata in $X$ are locally convex Riemannian
manifolds, and their connected components can be recognized metrically as
points in $X$ with isometric tangent spaces.  
The maximal dimensional stratum $X_{reg}$ is
unique, convex, open and dense, and consists of the \emph{principal orbits},
namely, those orbits with minimal isotropy groups, called
\emph{principal isotropy groups}. 
The \emph{cohomogeneity} of $\rho$ is the codimension of the principal orbits,
which also coincides with the topological dimension of $X$. 
The closure of the
union of strata of codimension one in $X$ is called the \emph{boundary} of~$X$,
and it is denoted by~$\partial X$. We say that $p\in V$ is
\emph{$G$-principal} if it lies in a $G$-principal orbit, and it is
\emph{$G$-important} if it projects to a stratum of codimension one in $X$.
Locally at~$p\in V$, the orbit decomposition of $V$ is completely
determined by the~\emph{slice representation} of the isotropy group
$G_p$ on the normal space $\nu_p(Gp)$ to the orbit $Gp$: namely,
the tangent cone $T_xX$ is isometric to the orbit space of the slice
representation at~$p$, where $x$ is the projection of~$p$;
in particular, the fixed point
set of the slice representation is tangent to the stratum of $p$
in $V$, and the cohomogeneity of the slice representation modulo the fixed
point set is equal to the codimension of the stratum of $x$ in $X$.

The strategy to investigate reductions of representations introduced
in~\cite{GL} is based on the following dichotomy. Suppose
$\rho:G\to \OG V$ is a minimal reduction
of $\tau:H\to \OG W$, where $H$ is connected.
Then the principal isotropy groups of $\rho$ are trivial, and:
(i) either $G$ is connected, and $V/G=W/H$ has non-empty
boundary~\cite[Proposition~5.2]{GL};
(ii) or $G$ is disconnected, and $G/G^0$ is generated
by \emph{nice involutions},
that is, elements of $G$ of order $2$ that fix
a $G^0$-principal and $G$-important point, and act on $V/G^0$ as a
reflection~\cite[Proposition~1.2 and~\S~4.3]{GL}.

In case~(i),
the boundary components arise as $S^a$-isotropy types, where $a=1$ or $a=3$.
More precisely, if $p\in V$ is a $G$-important point, then 
$G_p$ is a sphere $S^a$ and the slice
  representation, modulo the fixed point set,  is $(S^1,\C)$ or $(S^3,\Q)$.
  In case~(ii), the boundary components arise as $S^0$-isotropy types, where
  $S^0=\Z_2$, and the slice representation at a corresponding important
  point, modulo the fixed point set,
  is $(\Z_2,\R)$. In any case, we have the dimension formula
  \begin{equation}\label{dim-form}
    \dim V- a - 1 = \dim G - \dim N_{G}(G_p) + \dim V^{G_p},
    \end{equation}
  where  $N_{G}(G_p)$ denotes the normalizer of~$G_p$ in $G$, and
  $\dim V^{G_p}$ denotes the fixed point set of $G_p$ in $V$. 

  \section{Abstract copolarity $7$}

  We proceed to prove Theorem~\ref{main} using the method outlined in
  section~\ref{prelim}. In this section, we address the case of
  abstract copolarity $7$. So assume $\tau:H\to \OG W$ is a non-polar
  irreducible a representation of a compact connected Lie group $H$
  that reduces to $\rho:G\to \OG V$, where $\dim G=7$ and $G^0$ is not
  Abelian, and $\rho$ is a minimal reduction. 

  Recall that $\rho$ is irreducible, since $\tau$ is assumed
  irreducible, and the irreducibility of a representation can be read
  off its quotient~\cite[\S~5.2]{GL}. If the restriction
  of $\rho$ to the identity component $G^0$ were reducible, then
  $\tau$ would be toric~\cite[Theorem~1.7]{GL}, so we know that also
  $\rho|_{G^0}$ is irreducible.
  Now the center of $G^0$ can be at most one-dimensional,
  so $G^0$ is locally isomorphic to $\U1\times\SU2\times\SU2$. 
  
  \subsection{Connected case}\label{7-c}
  Herein we assume $G$ is connected. Then $V/G$ has non-empty boundary.
  Let $p\in V$ be a $G$-important point. We first rule out the case
  $G_p=S^3$. Indeed a $S^3=\SU2$-subgroup of $G$ must correspond to
  one of the $\SU2$-factors in $G$, or the diagonal in $\SU2\times\SU2$. In any
  case, the nontrivial element $z$ in the center of $G_p$ is a central
  element of $G$. The fixed point set of $z$ is not zero, as it contains~$p$,
  and it is preserved by $G$; but this contradicts the irreducibility of~$G$.

  Now $G_p=S^1$. Since the rank of $G$ is $3$, $\dim N_G(G_p)\geq3$.
  By the dimension formula~(\ref{dim-form}),
  \[ \dim V \leq 6 + \dim V^{G_p}. \]
  Since $G$ contains a central circle, $\rho$ is a representation of complex
  type. For a maximal torus of $G$ containing $G_p$, we have that
  $V^{G_p}$ is a sum of weight spaces; more precisely,  
  those weight spaces whose weights lie in the annihilator of a generator
  of the Lie algebra of~$G_p$. Since the dimension of the
  annihilator is two, there are at most two linearly independent weights
  appearing in $V^{G_p}$. Next, note that the restriction of any weight
  to the central circle in $G$ is independent of the weight. It follows
  that there are no linearly dependent weights appearing in $V^{G_p}$.
  We deduce that there are at most two weights appearing in $V^{G_p}$; 
  hence
  \[ \dim V\leq10. \]
  Now the cohomogeneity of $\rho$ is at most $3$. It is known that
  such representations must be either polar or of abstract
  copolarity~$1$ (\cite[Theorem~5.1]{Str}, \cite[Theorem~1.1]{GOT}
  and~\cite[Corollary~1.6 and Example~1.9]{GL}). 
  Therefore this case is not possible. 

  \subsection{Disconnected case}\label{copol7-discon}
  Now we assume $G$ is disconnected. Then there exists a nice
  involution $w\in G$, namely, $w^2=I$ and
  \begin{equation}\label{dimform-discon}
    \dim V = 8 - \dim Z_{G^0}(w) + \dim V^w,
    \end{equation}
  where $Z_{G^0}(w)$ denotes the centralizer of $w$ in $G^0$ and
  $V^w$ denotes the fixed point set of~$w$ in $V$.

  The conjugation of $G^0$ by $w$ defines an involutive automorphism
  $\varphi$ of $G^0$, namely,
  \begin{equation}\label{varphi}
    \rho(\varphi(g)) = w\rho(g)w^{-1}
    \end{equation}
    for all $g\in G^0$. Since $w\not\in\rho(G^0)$, and the centralizer
    of $\rho(G^0)$ in $\OG V$ coincides with its center, $\varphi$ cannot be
    an inner automorphism of $G^0$.  


    We may lift $\rho|_{G^0}$ to the universal covering and assume
    $G^0=\U1\times\SU2\times\SU2$. Any automorphism of this group
    must preserve the central circle
  and permute the $\SU2$-factors. Since $\SU2$ has no
  outer automorphisms, the group of outer automorphisms
  $\mathrm{Aut}(G^0)/\mathrm{Inn}(G^0)$
  is generated by the commuting
  involutions $\sigma_1$, $\sigma_2:G^0\to G^0$, where
  $\sigma_1(t,g_1,g_2)=(\bar t,g_1,g_2)$ and $
  \sigma_2(t,g_1,g_2)=(z,g_2,g_1)$.
  We may replace $\sigma_1$
  by an appropriate $G^0$-conjugate and assume
  $\sigma_1(t,g_1,g_2)=(\bar t,\bar g_1,\bar g_2)$.

  Note that $V=\C\otimes_{\mathbb C}\C^m\otimes_{\mathbb C}\C^n$
as a $\U1\times\SU2\times\SU2$-representation.  
  We have $\sigma_1$ is induced by conjugation by $\epsilon\in \OG V$,
  where $\epsilon$ is complex conjugation of $V$ over the real form
$\R\otimes_{\mathbb R}\R^m\otimes_{\mathbb R}\R^n$, namely,
  \begin{equation}\label{epsilon}
    \rho(\sigma_1(g))=\epsilon\rho(g)\epsilon^{-1}
    \end{equation}
  for all $g\in G^0$. 

In case $m=n$, 
also $\sigma_2$ is induced by conjugation by $\iota\in\OG V$, where
$\iota(z,z_1,z_2)=(z,z_2,z_1)$, namely,  
\begin{equation}\label{iota}
  \rho(\sigma_2(g))=\iota\rho(g)\iota^{-1}
  \end{equation}
  for all $g\in G^0$.

  Write $\varphi=\mathrm{Inn}_h\circ\sigma$, where
  $1\neq\sigma\in\langle\sigma_1,\sigma_2\rangle$ and $h\in G^0$,
  and $\mathrm{Inn}_h$ denotes the inner automorphism defined by~$h$. 
  We may assume $h$ has no component in the central circle. 
  Now $\varphi^2=I$ implies $h\sigma(h)$ lies in the center
  $Z(G^0)$ of $G^0$. Since either $(1,-1,1)$ or $(-1,-1,1)$ 
(resp.~$(1,1,-1)$ or $(-1,1,-1)$)
  lies in the
  kernel of~$\rho$, we can write
  \begin{equation}\label{h}
    \sigma(h)=h^{-1}.
  \end{equation}

  If $\varphi$ permutes the $\SU2$-factors of $G^0$, i.~e.~$\sigma$
  involves $\sigma_2$, then $m=n$. In particular, if 
  $\sigma=\sigma_2$, then~(\ref{h}) implies that
  $h=(1,h_1,h_1^{-1})\in G_0$ for some $h_1\in\SU2$. Replacing $w$ by
  $\rho(1,1,h_1)w\rho(1,1,h_1)^{-1}$,
  we may assume that $\varphi=\sigma_2$. It follows from~(\ref{varphi})
  and~(\ref{iota}) that
  $w\iota$ centralizes $\rho(G_0)$. By Schur's lemma,
  $w=\lambda\iota$ for some $\lambda\in\U1$, and from $w^2=I$ we deduce that
  $w=\pm\iota$. In any case, $\dim Z_{G^0}(w)=4$, so~(\ref{dimform-discon})
  gives $2n^2=4+\dim V^w$. Note that $V^w$ is the space of symmetric (resp.
  skew-symmetric) tensors, so $\dim V^w=n(n+1)$ (resp. $n(n-1)$). It follows
  that $n^2\pm n=4$, but these equations have no integer solutions, so 
  $\sigma=\sigma_2$ is impossible. 

  Similarly,  if $\sigma=\sigma_1\sigma_2$, then~(\ref{h})
implies that $h=(1,h_1,\bar h_1^{-1})\in G_0$ for some $h_1\in\SU2$. Replacing $w$ by
$\rho(1,1,\bar h_1)w\rho(1,1,\bar h_1)^{-1}$,
we may assume that $\varphi=\sigma_1\sigma_2$. It follows from~(\ref{epsilon})
  and~(\ref{iota}), together with Schur's lemma and $w^2=I$,
  that $w=\pm\iota\epsilon$. In any case, $\dim Z_{G^0}(w)=3$ so~(\ref{dimform-discon})
  gives $2n^2=5+\dim V^w$. Note that $V^w$ is the space of Hermitean (resp.
  skew-Hermitean) tensors, so $\dim V^w=n^2$. It follows
  that $n^2=5$, but this equation has no integer solutions, so 
  $\sigma=\sigma_1\sigma_2$ is impossible. 

  It remains the case $\sigma=\sigma_1$. Using~(\ref{h}), we obtain
  $h=(1,h_1,h_2)$, where $h_i\in\SU2$ is a symmetric matrix, for $i=1$, $2$.
  Taking the imaginary part of $h_ih_i^*=h_i\bar h_i=1$ shows that
  $\Re h_i$ and $\Im h_i$ are
  commuting (real symmetric) matrices. Therefore there is $k_i\in\SO2$ such that
  $k_ih_ik_i^{-1}$ is a diagonal matrix $d_i$. Now let $d_i^{1/2}$ be a
  square root of $d_i$ and put $a_i=d_i^{-1/2}k_i$. Then $a_ih_ia_i^t=1$.
  Let $a=(1,a_1,a_2)$. Then $\varphi(a)=ha$. Using this, we see that 
  by replacing $w$ by $\rho(a)w\rho(a)^{-1}$,
  $\varphi$ gets replaced by $\sigma_1$.
  Therefore we may assume $\varphi=\sigma_1$. Now $w\epsilon$ centralizes
  $\rho(G^0)$, so Schur's lemma and $w^2=I$ yield that $w=\pm\epsilon$.
  In any case, $\dim Z_{G^0}(w)=2$ and $\dim V^w=\frac12\dim V$,
  so~(\ref{dimform-discon}) yields that $\dim V=12$. This is
  the representation $(\U2\times\SU2,\C^2\otimes_{\mathbb C}\C^3)$ or,
  equivalently, $(\U2\times\SO3,\R^4\otimes_{\mathbb R}\R^3)$,
  which is the representation in the statement of Theorem~\ref{main}.
  
  \section{Abstract copolarity $8$}
 In this section we assume $\tau:H\to \OG W$ is a non-polar
  irreducible representation of a compact connected Lie group $H$
  that reduces to $\rho:G\to \OG V$, where $\dim G=8$ and $G^0$ is not
  Abelian, and $\rho$ is a minimal reduction. As in the previous
  section, $\rho|_{G^0}$ must be irreducible, and hence $G^0$ must be
  covered by $\SU3$.

  \subsection{Connected case}\label{8-c}
  Herein we assume $G$ is connected. We will show that $V/G$ must have empty
  boundary and hence $\rho$ cannot be a reduction of $\tau$. This follows
  immediately from the following result from~\cite{GKW}: \emph{Every
    irreducible representation of a compact connected simple Lie
    group with non-empty boundary in the orbit-space must be polar,
    toric, q-toric or a half-spin representation of $\Spin{11}$}.
  In the sequel we give an alternative argument. 

  We first note that there are no $S^3$-boundary components.
  Indeed it follows follows
  from~\cite[Corollary~13.4]{S} that a necessary condition
  for the existence of boundary is that the Dynkin index of the
  complexification $\rho^c$ is less than one.
  But~\cite[Table~1]{A-E-V} implies that such a representation
  would be polar. 

  Now, let $p$ be a $G$-important point projecting to a $S^1$-boundary
  component. Since the rank of $G^0$ is $2$, the dimension
  formula~(\ref{dim-form}) says
  \begin{equation}\label{copol8-s1}
    \dim V\leq 8 + \dim V^{G_p}.
    \end{equation}
  We next estimate $\dim V^{G_p}$. Consider a maximal torus $T$ of $G^0$ that
  contains $G_p$, so that $V^{G_p}$ is a sum of weight spaces. Recall
  that the diagram of weights of $V$ look like a sequence of concentric hexagons
  followed by a sequence of concentric triangles~\cite[p.~183-184]{FH}.
  Since
  the Lie algebra of $G_p$ has codimension one in the Lie algebra of $T$,
  the weights associated to weight spaces
  appearing in~$V^{G_p}$ all lie in a line through the origin. This line
  can meet at most two weights in each hexagon (resp.~triangle).
  
  Consider first the case $\rho$ is the realification of $\pi_{a,b}$,
  where $a>b$,
  $\pi_{a,b}$ is the complex irreducible representation of $\SU3$
  of highest weight $a\lambda_1+b\lambda_2$, and $\lambda_1$, $\lambda_2$
  denote the fundamental weights. Then the diagram of weights
  consists of hexagons $H_0$,\ldots, $H_{b-1}$, and triangles
$T_0,\ldots,T_{[\frac{a-b}3]}$, the multiplicities of
  weights on $H_i$ being $i+1$, and the multiplicities of weights
  on any triangle being $b+1$. We obtain the rough estimate
  \[ \dim V^{G_p} \leq 2b(b+1) +\frac43(b+1)(a-b). \] 
  The Weyl dimension formula gives $\dim V=(a+1)(b+1)(a+b+2)$,
  and one easily sees that the only solution to~(\ref{copol8-s1}) is
  $(a,b)=(1,0)$, which gives the vector representation,
  a polar representation.

  Consider next the case $\rho$ is a real form of~$\pi_{a,a}$. Then
  the diagram of weights consists of hexagons $H_0,\ldots,H_{a-1}$,
  the multiplicities of  weights on $H_i$ being $i+1$,
  and a point $H_a$ with multiplicity $a+1$. We obtain
  \[ \dim V^{G_p} \leq (a+1)^2. \]
  Since $\dim V=(a+1)^3$, the only solution to~(\ref{copol8-s1})
  is $a=1$, the adjoint representation, which is polar. 

  \subsection{Disconnected case}
Herein we assume $G$ is disconnected. We will show that $V/G$ must have empty
boundary and hence $\rho$ cannot be a reduction of $\tau$.

There is a nice involution $w\in G$, so that $w^2=I$ and
\begin{equation}\label{df-8-d}
  \dim V = 9 - \dim Z_{G^0}(w) + \dim V^w.
  \end{equation}
As in subsection~\ref{copol7-discon}, the conjugation of $G^0$ by~$w$
defines an automorphism $\varphi$ of $G^0$. We will deal
separately with the cases in which $\varphi$ is of inner or outer type. 

\subsubsection{Outer type}\label{8-d-o}
It follows from~\cite{W} that, replacing $w$ by a $\rho(G^0)$-conjugate,
we may assume $\varphi=\sigma$, where $\sigma$ is the Weyl involution,
given by $\sigma(g)=\bar g$ for all $g\in G^0$.

We first check $\rho|_{G^0}$ cannot be absolutely irreducible. Indeed, if this
is the case, the complexification $(\rho|_{G_0})^c=:\pi$ is irreducible. 
Denote by $\epsilon$ the complex conjugation of $V^c$ over $V$. Then
$\epsilon\circ\pi\circ\epsilon$ and $\pi\circ\sigma$ are equivalent
representations. So there is a unitary transformation $A$ of $V^c$
such that
\[ A\epsilon\pi(g)\epsilon A^{-1}=\pi(\sigma(g))=w\pi(g)w^{-1}, \]
for all $g\in G^0$, where we have denoted the complex linear extension
of $w$ to $V^c$ by the same letter. This says that $(A\epsilon)^{-1}w$
centralizes $\pi(G^0)$, so Schur's lemma yields that $w$ is a complex
multiple of $A\epsilon$. We reach a contradiction, as $\epsilon A$ is not
a complex transformation of $V^c$.

Now  $\rho|_{G^0}$ is not absolutely irreducible. Then it is
the realification of a representation $\pi$ of $G^0$ on a complex
vector space~$U$. Let $\epsilon$ be the
conjugation of $U$ over any real form.  Then $\epsilon\circ\pi\circ\epsilon$ and $\pi\circ\sigma$ are equivalent
representations, so as above we can write $w=\lambda\epsilon A$, for
some $\lambda\in S^1$, and some unitary transformation~$A$ of~$U$. 
We can replace $w$ by $\lambda^{-1/2}w\lambda^{1/2}=\epsilon A$.  
Next, since $\epsilon$ is conjugate linear, we have
\[ I = w^2 =(\lambda\epsilon A)^2 =(\epsilon A\epsilon)A = \bar A A, \]
that is, $A$ is unitary and symmetric. Write $A=BB^t$, where $B$ is a
unitary matrix, similar what was done in the end of
section~\ref{copol7-discon}. Then $B^tw(B^t)^{-1}=\bar B^{-1}\epsilon A(B^t)^{-1}
=\epsilon B^{-1}(BB^t)(B^t)^{-1}=\epsilon$, so we may indeed
assume $w=\epsilon$.
Since $w$ is conjugate linear, we have $\dim V^w=\frac12\dim V$.
Further, $\dim Z_{G^0}(w)=\dim \SO3=3$, so the dimension formula~(\ref{df-8-d})
yields $\dim V=12$. This implies that the cohomogeneity of~$\rho$, and hence
also of $\tau$, is~$4$,
but all non-polar, non-reduced, irreducible representations of
cohomogeneity~$4$ of compact connected Lie groups have abstract
copolarity~$2$~\cite[Theorem~1.11]{GL}. 

\subsubsection{Inner type}
Write $w=\rho(h)z$, where $h\in G^0$, and $I\neq z\in \OG V$ centralizes
$\rho(G^0)$. From $w^2=I$ we deduce that $h^2$ lies in the center $Z(G^0)$.
If $h\in Z(G^0)$ then $\rho(h)$ and $z$ are scalar maps, and so $w=\pm I$,
which cannot be. Now $h\not\in Z(G^0)$, and we can
conjugate $h$ to a diagonal matrix and assume
$h=\left(\begin{smallmatrix}\theta&0&0\\0&-\theta&0\\0&0&-\theta\end{smallmatrix}\right)$, where
$\theta$ is a cubic root of~$1$. It follows that
$\dim Z_{G^0}(w)=\dim Z_{G^0}(h)=4$, and the dimension formula gives $\dim V-\dim V^w=5$
is odd. This already implies that $\rho$ cannot leave invariant 
a complex structure on $V$.
Now Schur's lemma yields $z=-I$, so 
$w=-\rho(h)=-\rho(\theta\left(\begin{smallmatrix}1&0&0\\0&-1&0\\0&0&-1\end{smallmatrix}\right))=\pm\rho(\left(\begin{smallmatrix}1&0&0\\0&-1&0\\0&0&-1\end{smallmatrix}\right))$, 
and we get
$\dim V^h=5$.

We know that $\rho$ is a real form of $\pi_{a,a}$ (notation as
in subsection~\ref{8-c}). Consider the maximal torus of $G^0$
consisiting of diagonal matrices, so $(V^c)^h=(V^h)^c$
is a sum of weight spaces. It is clear that the set of weights
associated to weight spaces in  $(V^c)^h$ is invariant under
multiplication by $-1$. Since the multiplicity of the zero weight
in $\pi_{a,a}$ is $a+1$, we deduce that
$5=\dim V^h=\dim_{\mathbb C}(V^c)^h\geq a+1$, and~$a$ cannot be an odd number.
We are only left to examine the cases $a=2$ and $a=4$.

Let $e_1$, $e_2$, $e_3$ denote the canonical basis of $\C^3$
and let $e_1'$, $e_2'$, $e_3'$ denote the dual basis of $\C^{3*}$.
The highest weight vector of $\pi_{a,a}$ is $e_1^a\otimes e_1'^a$,
which is fixed by $h$ for all $a$. Using the
action of the Weyl group, the vectors $e_j^a\otimes e_j'^a$
lie in $V^c$ for $j=1$, $2$, $3$, and they are clearly also
fixed by $h$, together with the corresponding weight
vectors with the opposite weights, and the zero weight space.
This gives the estimate
$\dim V^h\geq 3\cdot 2 + (2+1)=9$
for $a\geq2$, which contradicts $\dim V^h=5$. 

\section{Abstract copolarity~$9$}
In this section we assume $\tau:H\to \OG W$ is a non-polar
  irreducible representation of a compact connected Lie group $H$
  that reduces to $\rho:G\to \OG V$, where $\dim G=9$ and $G^0$ is not
  Abelian, and $\rho$ is a minimal reduction. As in the previous
  section, $\rho|_{G^0}$ must be irreducible, and hence $G^0$ must be
  locally isomorphic to~$\U3$ or $\SP1^3$, but the latter case is q-toric, so
  we need not consider it. Now $\rho$ is a representation of complex type.  

 \subsection{Connected case}\label{9-c}
 Herein we assume $G$ is connected, so $G$ is covered by $\U1\times\SU3$.
 We will show that $V/G$ must have empty
  boundary and hence $\rho$ cannot be a reduction of~$\tau$.

  Note that if $p$ is a $G$-important point that projects to
  a $S^3$-boundary component, then the restriction of $\rho$
  to the simple subgroup~$\SU3$
  is a (non-orbit equivalent)
  representation which is either irreducible of complex type,
  or with two equivalent irreducible components of real type,
  and which also has~$p$ as an important
  point projecting to a $S^3$-boundary component.
  By~\cite[\S13]{S}
  and~\cite[Table~1]{A-E-V}, there are no possibililties for~$\rho$.

  Now, let $p$ be a $G$-important point projecting to a $S^1$-boundary
  component. Since the rank of $G$ is $3$, the dimension
  formula~(\ref{dim-form}) says
  \begin{equation}\label{copol9-s1}
    \dim V\leq 8 + \dim V^{G_p}.
  \end{equation}
  Again, since $G$ has rank $3$ and contains a central circle, as
  in subsection~\ref{7-c}
    we see that $V^{G_p}$ is a sum of at most two weight spaces.
    Suppose $\rho$ is the realification of $\pi_{a,b}$. Then the
    highest multiplicity of a weight is~$b+1$, so $\dim V^{G_p}\leq 4b+4$,
    and~(\ref{copol9-s1}) implies
    \[ (a+1)(b+1)(a+b+2) \leq 12 + 4b. \]
    The only solution is~$(a,b)=(1,0)$, which corresponds to a polar
    representation.

    \subsection{Disconnected case}
    Herein we assume $G$ is disconnected.
        We will show that $V/G$ must have empty
boundary and hence $\rho$ cannot be a reduction of $\tau$.

There is a nice involution $w\in G$, so that $w^2=I$ and
\begin{equation}\label{df-9-d}
  \dim V = 10 - \dim Z_{G^0}(w) + \dim V^w.
  \end{equation}
As in subsection~\ref{copol7-discon}, the conjugation of $G^0$ by~$w$
defines an automorphism $\varphi$ of $G_0$. Since $w\not\in\rho(G^0)$,
and the centralizer
    of $\rho(G^0)$ in $\OG V$ coincides with its center, $\varphi$ must be
    an outer automorphism of $G^0$.

    Since the center of $\rho(G^0)$ acts by scalars and is
    normalized by~$w$, 
    we may assume that $G^0=\U3$, up to orbit-equivalence.
    Now the group of outer automorphisms of $G^0$ is generated
    by $\sigma$, where $\sigma(g)=\bar g$ for $g\in\U3$. 
    Similar to the situation in subsection~\ref{8-d-o}, we may replace
$w$ by a conjugate and assume $\varphi=\sigma$.
Here we already know $\rho|_{G^0}$ is not absolutely irreducible. 
Continuing with the argument in  subsection~\ref{8-d-o},
we may assume $w=\epsilon$. Now $\dim Z_{G^0}(w)=3$,
so the dimension formula~(\ref{df-9-d})
yields $\dim V=14$. However, there exist no irreducible
representations of dimension~$14$ of~$\U3$. This finishes the proof
of Theorem~\ref{main}.


\begin{thebibliography}{GKW21}

\bibitem[AV{\'E}67]{A-E-V}
E.~M. Andreev, E.~B. Vinberg, and A.~G. {\'E}lashvili, \emph{Orbits of greatest
  dimension in semi-simple linear {L}ie groups. ({R}ussian)}, {Funktsional.\
  Anal.\ i Prilozhen.\/} \textbf{1} (1967), no.~4, 3--7, English transl.,
  Functional Analysis and Its Applications, 1967, 1:4, 257--261.

\bibitem[BCO03]{BCO}
J.~Berndt, S.~Console, and C.~Olmos, \emph{Submanifolds and holonomy}, Research
  Notes in Mathematics, no. 434, {Chapman \& Hall/CRC, Boca Raton}, 2003.

\bibitem[Dad85]{D}
J.~Dadok, \emph{Polar coordinates induced by actions of compact {L}ie groups},
  {Trans. Amer. Math. Soc.} \textbf{288} (1985), 125--137.

\bibitem[DK85]{DK}
J.~Dadok and V.~Kac, \emph{Polar representations}, {J.\ Algebra} \textbf{92}
  (1985), 504--524.

\bibitem[FH91]{FH}
W.~Fulton and J.~Harris, \emph{Representation theory. {A} first course},
  Graduate Texts in Mathematics, vol. 129, Springer-Verlag, New York, 1991.

\bibitem[GG18]{GG}
C.~Gorodski and F.~Gozzi, \emph{{Representations with $Sp(1)^k$-reductions and
  quaternion-K\"ahler symmetric spaces}}, {Math.\ Z.\/} \textbf{290} (2018),
  561--575.

\bibitem[GKW21]{GKW}
C.~Gorodski, A.~Kollross, and B.~Wilking, \emph{Action on positively curved
  manifolds and boundary in the orbit space}, work in progress, 2021.

\bibitem[GL14]{GL}
C.~Gorodski and A.~Lytchak, \emph{On orbit spaces of representations of compact
  {L}ie groups}, {J.\ reine angew.\ Math.\/} \textbf{691} (2014), 61--100.

\bibitem[GL15]{GL2}
\bysame, \emph{Representations whose minimal reduction has a toric identity
  component}, {Proc.\ Amer.\ Math.\ Soc.\/} \textbf{143} (2015), 379--386.

\bibitem[Gom21]{Gom}
A.~M. de~S\'a Gomes, \emph{Representations of low copolarity and biquotients},
  Ph.D. thesis, University of S{\~a}o Paulo, 2021.

\bibitem[GOT04]{GOT}
C.~Gorodski, C.~Olmos, and R.~Tojeiro, \emph{Copolarity of isometric actions},
  {Trans.~Amer.~Math.~Soc.} \textbf{356} (2004), 1585--1608.

\bibitem[Pan17]{Pa}
F.~Panelli, \emph{Representations admitting a toric reduction}, Ph.D. thesis,
  University of Florence, 2017.

\bibitem[Sch80]{S}
G.~W. Schwarz, \emph{Lifting smooth homotopies of orbit spaces}, {I.H.E.S.
  Publ. in Math.} \textbf{51} (1980), 37--135.

\bibitem[Str94]{Str}
E.~Straume, \emph{On the invariant theory and geometry of compact linear groups
  of cohomogeneity $\leq3$}, {Diff. Geom. and its Appl.} \textbf{4} (1994),
  1--23.

\bibitem[Wol11]{W}
J.~A. Wolf, \emph{Spaces of constant curvature}, sixth ed., AMS Chelsea
  Publishing, Providence, RI, 2011.

\end{thebibliography}

\providecommand{\bysame}{\leavevmode\hbox to3em{\hrulefill}\thinspace}
\providecommand{\MR}{\relax\ifhmode\unskip\space\fi MR }
\providecommand{\MRhref}[2]{%
  \href{http://www.ams.org/mathscinet-getitem?mr=#1}{#2}
}
\providecommand{\href}[2]{#2}

\end{document}